# An ancient Egyptian problem:
# The diophantine equation

$$\frac{4}{n} = \frac{1}{x} + \frac{1}{y} + \frac{1}{z}, n \geq 2$$

*Konstantine Zelator*
*Department of Mathematics*
*301 Thackeray Hall*
*139 University Place*
*University of Pittsburgh*
*Pittsburgh, PA 15260*
*U.S.A*

*Also: Konstantine Zelator*
*P.O. Box 4280*
*Pittsburgh, PA 15203*
*U.S.A*

*e-mail addresses: 1) konstantine_zelator@yahoo.com*
*2) kzet159@pitt.edu*

INTRODUCTION

From the Rhind Papyrus and other sextant sources, we know that the ancient Egyptians were very interested in expressing a given fraction into a sum of distinct unit fractions, that is, fractions whose numerators are equal to 1. They even developed tables where in they tabulated the unit fraction decompositions of specific fractions. One of the problems that has come down to us in the last 60 years or so; is the problem of whether for each $n \geq 2$; the fraction $\frac{4}{n}$ can be decomposed as a sum of three distinct unit fractions. In otherwords, whether for each $n \geq 2$; the diophantine equation $\frac{4}{n} = \frac{1}{x} + \frac{1}{y} + \frac{1}{z}$ has a solution in positive integers $x, y$ and $z$; with $x \neq y$, $y \neq z$ and $z \neq x$. This problem is formally known as the Erdos-Strauss Conjecture, first formulated in 1948, even though the earliest published reference to it appears to be a 1950 Paul Erdos paper. Since 1950, a number of partial results have been obtained, for example see references **[1] – [10], [12]** and **[13].** In this work, we contribute four theorems, three of which (Theorems 2,3, and 4) directly deal with the above problem.

## 2. THE THEOREMS

**THEOREM 1:**

(a) Let $p$ and $q$ be natural numbers, with $p$ a prime. Then, the diophantine equation $\frac{q}{p} = \frac{1}{x_1} + \frac{1}{x_2}$ has a solution in the set of natural numbers $N$, with $x_1$ and $x_2$ distinct, if and only if $p + 1 \equiv 0 \pmod{q}$. If the above equation is solvable, then the solution is uniquely determined (up to symmetry) and given by $x_1 = \frac{p+1}{q}$ and $x_2 = p\left(\frac{p+1}{q}\right)$.

(b) If $p$ and $q$ are arbitrary positive integers with $p+1 \equiv 0 (\mod q)$, then $\left\{ \dfrac{p+1}{q}, p\left(\dfrac{p+1}{q}\right) \right\}$ is a solution to the same equation.

**THEOREM 2:**

If the diophantine equation $\dfrac{4}{n} = \dfrac{1}{x_1} + \dfrac{1}{x_2} + \dfrac{1}{x_3}$ has no solution in $N$ with $x_1, x_2$ and $x_3$ distinct, then every prime divisor $p$ of $n$ must satisfy $p \equiv 1(\mod 24)$ $(and\ so\ n \equiv 1(\mod 24))$. Also, if $n$ has a divisor not congruent to $1(\mod 24)$, then a solution to the above equation can be explicitly found.

**THEOREM 3:**

Let $n$ be an odd natural number with no divisors $\equiv 3(\mod 4)$, and $\delta$ a divisor of $n$. Suppose that there exists an odd natural number $k$ such that $\delta + k \equiv 0(\mod m)$, for some positive integer $m \equiv 3(\mod 4)$ and such that $n\left(\dfrac{m+1}{4}\right)\left(\dfrac{\delta+k}{m}\right) \equiv 0(\mod k)$. Then the diophantine equation $\dfrac{4}{n} = \dfrac{1}{x_1} + \dfrac{1}{x_2} + \dfrac{1}{x_3}$ has a solution in $N$, with $x_1, x_2$, and $x_3$ distinct.

**REMARK:** The above theorem can also be stated as follows: If has diophantine system $\delta + k = a(4t-1), atn \equiv 0(\mod k)$ and $n \equiv 0(\mod \delta)$ is solvable in $N$, then the above diophantine equation has a solution in $N$ (with $x_1, x_2$ and $x_3$ distinct; note that $t = \dfrac{m+1}{4}$).

**THEOREM 4:**

Let $n$ be an odd natural number with no divisors $\equiv 3(\mod 4)$, and suppose that there exist divisors $\delta$ and $d$ of $n$ such that $\delta + d$ contains a divisor congruent to $3$ modulo $4$. Then the diophantine equation $\dfrac{4}{n} = \dfrac{1}{x_1} + \dfrac{1}{x_2} + \dfrac{1}{x_3}$ has a solution in $N$, with $x_1, x_2$ and $x_3$ distinct. If in particular $n$ has a divisor $\delta$ such that $\delta + 1 \equiv 0(\mod m)$, for some $m \equiv 3(\mod 4)$, then such a solution exists.

## 3. THE PROOFS

**PROOF OF THEOREM 1:**

(a) Let $\{x_1, x_2\}$ be a solution of distinct positive integers $x_1$ and $x_2$ to the equation,

$$\frac{q}{p} = \frac{1}{x_1} + \frac{1}{x_2} \qquad (1)$$

Writing (1) in the form $qx_1 x_2 = p(x_1 + x_2)$ \qquad (2),

we see that since $p$ is a prime, it must divide $q$, $x_1$, or $x_2$. It cannot divide $q$ for in such a case the fact that $\frac{q}{p} \geq 1$ would imply $\frac{q}{p} x_1 x_2 \geq x_1 x_2 \geq x_1 + x_2$, since $x_1$ and $x_2$ are both $\geq 1$; if both $x_1$ and $x_2$ are greater than $1$, then $x_1 x_2 > x_1 + x_2$ and so the above inequality implies $\frac{q}{p} x_1 x_2 > x_1 + x_2$, in violation of (2); if say $x_1 = 1$, equation (2) implies $\frac{q}{p} x_2 = 1 + \frac{1}{x_2}$ which is a contradiction since the left-hand side of the last equation is an integer but the right-hand side cannot be, because $x_2 > 1$ ($x_1$ and $x_2$ are distinct).

So $p$ must divide $x_1$ or $x_2$. Without loss of generality let us assume that $p$ divides $x_1$ so that $\frac{x_1}{p}$ is a natural number. So from equation (2) $\Rightarrow q\left(\frac{x_1}{p}\right) x_2 - p\left(\frac{x_1}{p}\right) = x_2$ \qquad (3)

Equation (3) shows that the integer $\frac{x_1}{p}$ divides $x_2$, so let us set

$$x_2 = \left(\frac{x_1}{p}\right) k \qquad (4),$$

$k$ a pos.integer.

Combining (3) and (4) we obtain,

$$kq\left(\frac{x_1}{p}\right) = p + k \qquad (5)$$

Equation (5) gives

$$x_1 = \frac{p(p+k)}{kq} \tag{6}$$

Since $x_1$ is an integer, (6) implies that every prime divisor of $k$ must divide $p$ or $p+k$; but obviously then it follows that, since $p$ is a prime, we must have $k = p^\alpha$, where $\alpha$ is some pos.integer. Since $p^\alpha$ divides $p(p+k) = p(p+p^\alpha) = p^2(1+p^{\alpha-1})$, it is clear that $\alpha \leq 2$. For $\alpha = 1$ we have $k = p^\alpha = p$ and so eq. (6) yields $x_1 = \frac{2p}{q}$ which implies that, since $p$ is a prime and $x_1$ an integer, $q = 1, 2, p$ or $2p$. But we have already seen that $p$ cannot divide $q$ (refer to the earlier argument involving eq. (2)). Thus $q = 1$ or $2$. If $q = 1$, then trivially $p+1 \equiv 0 \pmod{q}$; and from (6) and (4) we have $x_1 = x_2 = 2p$ which violates the condition $x_1 \neq x_2$; if $q = 2$ then from $k = p$ and eq. (6) and (4) we obtain $x_1 = x_2 = p$, again violating $x_1 \neq x_2$. Now assume $\alpha = 2 \Rightarrow k = p^2$; and by eq. (6) we obtain,

$$x_1 = \frac{p+1}{q} \tag{7}$$

and

from eq. (4) $\Rightarrow$ $$x_2 = \left(\frac{p+1}{q}\right)p \tag{8}$$

Thus $p+1 \equiv 0 \pmod{q}$. The converse of the Theorem is immediate, for if $x_1$ and $x_2$ are integers given by (7) and (8), they clearly $x_1 \neq x_2$, and a straightforward computation shows that $\{x_1, x_2\}$ is a solution to equation (1).

(b) The proof is immediate by direct computation.

**Proof of Theorem 2:**

We will demonstrate that if $n$ is a positive integer with at least one prime divisor $p \neq 1 \pmod{24}$, then the equation in question has a solution with $x_1, x_2$ and $x_3$ distinct. If $n \equiv -1 \pmod 3$, then the equation

$$\frac{4}{n} = \frac{1}{x_1} + \frac{1}{x_2} + \frac{1}{x_3} \tag{9}$$

is indeed solvable with $x_3 = n$, since the equation $\dfrac{3}{n} = \dfrac{1}{x_1} + \dfrac{1}{x_2}$ has by Th. 1 part (b), the solution $x_1 = \dfrac{n+1}{3}, x_2 = \left(\dfrac{n+1}{3}\right)$; in the notation of Th.1, $q = 3$ and $p = n$. If $n \equiv 0 \pmod{3}$, then again take $x_3 = n$; equation (9) becomes

$$\dfrac{3}{n} = \dfrac{1}{x_1} + \dfrac{1}{x_2}; \dfrac{1}{n/3} = \dfrac{1}{x_1} + \dfrac{1}{x_2},$$

which by Th. 1 has the solution $x_1 = \dfrac{n}{3} + 1$ and $x_2 = \left(\dfrac{n}{3} + 1\right)\dfrac{n}{3}$ ( in the notation of Th.1, $p = n/3$ and $q = 1$).

We now come to the case $n \equiv 1 \pmod{3}$. To simplify matters, we observe that if $n$ is even, (9) has a solution: to see this take $x_3 = n/2$ and so (9) $\Leftrightarrow \dfrac{1}{n/2} = \dfrac{1}{x_1} + \dfrac{1}{x_2}$, which of course has the solution $x_1 = n/2 + 1$ and $x_2 = (n/2 + 1)n/2$. Thus we may assume that $n \equiv 1 \pmod{3}$ and $n \equiv 1 \pmod{2}$. In effect, $n \equiv 1 \pmod{6}$. First assume that $\dfrac{n-1}{6}$ is odd. We can write (9) in the form,

$$(4x_3 - n)x_1 x_2 = nx_3(x_1 + x_2) \qquad (10)$$

If $n$ were prime $n$ would have to divide at least one of the numbers $4x_3 - n, x_1$, and $x_2$; suppose $4x_3 - n \equiv 0 \pmod{n}$; then $x_3 \equiv 0 \pmod{n}$ and put $x_3 = \delta n$, $\delta$ some natural number. If in addition we set $x_2 = \delta x_2'$, then equation (10) would imply,

$$(4\delta - 1)x_1 x_2' = n(x_1 + \delta x_2') \qquad (11)$$

Moreover if we put $x_1 = x_2' = 1 + \delta$ in (11), we obtain

$$4\delta - 1 = n \text{ which is true for } \delta = \dfrac{n+1}{4},$$

in virtue of the fact that $\dfrac{n-1}{6}$ is odd. Consequently the numbers $x_1 = 1 + \delta = \dfrac{n+5}{4}, x_2 = \delta(1+\delta) = \dfrac{(n+1)(n+5)}{16}$, and $x_3 = \delta n = \dfrac{n(n+1)}{4}$, constitute a solution to (9), when $n \equiv 1 \pmod 6$ and $\frac{n-1}{6}$ is odd. Next suppose that $n \equiv 1 \pmod 6$ and $\frac{n-1}{6} \equiv 0 \pmod 2$; that is, assume $n \equiv 1 \pmod{12}$ and let $p$ be any prime divisor of $n$.

If $p \not\equiv 1 \pmod{12}$, then we have already shown that the equation $\frac{4}{p} = \frac{1}{x_1} + \frac{1}{x_2} + \frac{1}{x_3}$ has a solution in $N$ with $x_1, x_2,$ and $x_3$ distinct. By multiplying the last equation with $\frac{1}{n/p}$ we obtain a solution $x_1' = x_1 \, n/p, x_2' = x_2 \, n/p, x_3' = x_3 \, n/p$ to equation (9). If $p \equiv 1 \pmod{12}$ and $p \not\equiv 1 \pmod{24}$, then $\frac{p-1}{12}$ is an odd integer. The integer $p$ can be expressed in the form $p = 4k - 3$, with $k \equiv 4 \pmod 6$. We have,

$$\frac{4}{n} - \frac{1}{kn} = \frac{4k-1}{nk} = \frac{4k-3}{nk} + \frac{2}{nk} = \frac{1}{\left(\frac{n}{p}\right)\left(\frac{p+3}{4}\right)} + \frac{1}{n\left(\frac{p+3}{8}\right)},$$

since $k = \frac{p+3}{4}$. Therefore we conclude that

$$\frac{4}{n} = \frac{1}{\left(\frac{p+3}{4}\right)n} + \frac{1}{\left(\frac{n}{p}\right)\left(\frac{p+3}{4}\right)} + \frac{1}{n\left(\frac{p+3}{8}\right)}, \text{ which}$$

shows the equation (9) has a solution with distinct $x_1, x_2, x_3$, and the proof is concluded.

### PROOF OF THEOREM 3

According to hypothesis, there exists $m = 4t - 1$ and $k$ odd such that

$\delta + k = a(4t - 1)$ and $atn \equiv 0 \pmod k$, (12), for some number $a$; $\delta$ a divisor of $n$. Now, observe that

$$\frac{4}{n} - \frac{k}{atn} = \frac{4at - k}{atn}, \text{ and by}$$

applying the equation in (12) we obtain

$$\frac{4}{n} - \frac{k}{atn} = \frac{\delta + a}{atn} = \frac{\delta}{atn} + \frac{1}{tn}$$

$$\Rightarrow \frac{4}{n} = \frac{1}{atn/k} + \frac{1}{at(n/\delta)} + \frac{1}{tn} \qquad (13)$$

However, according to the congruence in (12), $atn/k$ is a natural number, and so is $n/\delta$, by virtue of $\delta$ being a divisor of $n$. Consequently (13) shows that the numbers $x_1 = \dfrac{atn}{k}, x_2 = at\left(\dfrac{n}{\delta}\right)$ and $x_3 = tn$, form a solution to $\dfrac{4}{n} = \dfrac{1}{x_1} + \dfrac{1}{x_2} + \dfrac{1}{x_3}$. To prove that $x_1 = \dfrac{atn}{k}, x_2 = at\left(\dfrac{n}{\delta}\right)$ and $x_3 = tn$ are distinct integers, is an easy matter: if $x_1 = x_2$, then we must have $k = \delta$ and so from the equation in (12) it follows that $2\delta = a(4t - 1)$; but $\delta$ can not contain divisors congruent to $3 \pmod 4$, since $\delta$ is a divisor of $n$ and $n$ has no divisors congruent to $3 \pmod 4$. Secondly, if $x_1 = x_3$, we have $a = k$ and by (12) $\Rightarrow \delta = k[4t - 2]$, which is not possible since $\delta$ is odd, being a divisor of $n$ which is odd. Finally if $x_2 = x_3 \Rightarrow a = \delta$ and from (12) $\Rightarrow k = \delta(4t - 2)$, in violation of the fact that $k$ is odd.

**PROOF OF THEOREM 4:**

The proof is immediate if one observes that the hypothesis of Theorem 4 satisfies the hypothesis of Theorem 3 with $k = d$. So Th. 4 is really an obvious corollary of Th. 3. Note that the congruence $n\left(\dfrac{m+1}{4}\right)\left(\dfrac{\delta + k}{m}\right) \equiv 0 \pmod k$ in Th. 3 is obviously satisfied, since $k = d$ is divisor of $n$.